\documentclass[11pt,a4paper]{article}

\usepackage[margin=2.7cm]{geometry}
\usepackage{graphicx}
\usepackage{multirow}
\usepackage{amsmath,amssymb,amsfonts}
\usepackage{amsthm}
\usepackage{mathrsfs}
\usepackage[title]{appendix}
\usepackage{xcolor}
\usepackage{textcomp}
\usepackage{manyfoot}
\usepackage{booktabs}
\usepackage{algorithm}
\usepackage{algorithmicx}
\usepackage{algpseudocode}
\usepackage{listings}
\usepackage[hidelinks]{hyperref}

\theoremstyle{plain}
\newtheorem{theorem}{Theorem}
\newtheorem{proposition}[theorem]{Proposition}
\newtheorem{lemma}[theorem]{Lemma}
\newtheorem{corollary}[theorem]{Corollary}

\theoremstyle{remark}

\newtheorem{remark}{Remark}

\theoremstyle{definition}
\newtheorem{definition}{Definition}

\title{Pathwise Forcing of Exact Magnetization Levels in Glauber--Ising Dynamics}
\author{Jean-Gabriel Attali\thanks{\texttt{jean-gabriel.attali@devinci.fr}}\\
  De Vinci Higher Education, De Vinci Research Center\\
  Paris, France}
\date{}

\begin{document}
\maketitle

\begin{abstract}

For $d\geq2$, we show that every exact magnetization level
$C_r=\{\sigma:M_n(\sigma)\to r\}$, $r\in[-1,1]$,
of the heat-bath Glauber--Ising dynamics contains a configuration whose
empirical magnetization almost surely leaves $C_r$ at every
sufficiently small positive time. For $r\neq0$, forcing toward zero
occurs, while a nontrivial interval around zero admits forcing in both
directions at every positive temperature. Rational levels have
explicit periodic witnesses, and interpolation yields full-measure
families at arbitrary nonzero interior levels. Above a dimension-dependent
threshold, the outward direction is realized by biased Bernoulli
quenches. In dimension two, these quenches already exhibit nonmonotone
magnetization relaxation in part of the uniqueness regime. Throughout
the coexistence regime, the resulting two-way forcing interval contains
$[-m_\beta,m_\beta]$ and extends beyond it, producing a nonequilibrium
overshoot. In contrast, in dimension one $C_0$ is pathwise absorbing,
while every nonzero exact level is pathwise forced toward zero.

\end{abstract}

\section{Introduction}

We consider the nearest-neighbour Ising model on $\mathbb Z^d$,
$d\geq1$, at zero external field, under heat-bath Glauber dynamics~\cite{Glauber63}.
Our main results concern $d\geq2$, while the qualitatively different
one-dimensional case is treated separately.
General references for interacting particle systems,
graphical constructions, and attractive spin systems are
\cite{Harris72,Liggett85,Liggett99}; for Gibbs measures and the Ising
model, see \cite{Georgii11,FriedliVelenik17}. We use classical
two-dimensional equilibrium results from
\cite{Aizenman80,Higuchi81,Yang52}.

Nonequilibrium Ising dynamics has been studied in connection with
phase ordering and domain coarsening \cite{Bray1994}. 
At zero temperature, Fontes, Schonmann and Sidoravicius proved fixation to the all-plus state, with stretched-exponential relaxation bounds, from sufficiently biased product initializations 
\cite{FontesSchonmannSidoravicius02}. In dimension two, a
positive-temperature phase-ordering result was recently obtained
throughout the low-temperature regime \cite{GheissariSly26}. 
Dependence on the initial state also appears in finite-volume dynamics \cite{gheissariSinclair23} while information-percolation methods provide a different approach to stochastic Ising dynamics \cite{lubetzkySly16}. The evolution of Gibbsian structure under stochastic dynamics was studied in \cite{LeNyRedig02,vanEnterFernandezdenHollanderRedig02}.

Let $C_r$ denote the set of configurations whose empirical
magnetization along centered boxes converges to $r$. The value of $r$
fixes the asymptotic density of plus and minus spins, but not the
frequencies of their local environments. The latter enter the
short-time empirical drift.


The one-dimensional model is highly special, at every finite
inverse temperature there is no phase coexistence, and the empirical
magnetization closes under the dynamics. The behavior of an entire exact level $C_r$ is then determined by $r$ alone. In sharp contrast, for $d\geq 2$ the same magnetization level can contain configurations with different local-environment frequencies and opposite short-time macroscopic drifts. We record the resulting complete one-dimensional picture for exact magnetization levels before turning to the higher-dimensional phenomenon.


\paragraph{Main result.}
For every $d\ge2$, every $\beta>0$, and every $r\in[-1,1]$, the level
$C_r$ contains a pathwise forcing witness: a deterministic
configuration whose trajectory almost surely lies strictly outside
$C_r$ at every sufficiently small positive time. For $r\ne0$, the
forcing direction can be chosen toward zero. At $r=0$, witnesses exist
in both directions. More generally, for every $\beta>0$, both
directions occur on a nontrivial interval around zero. Above a
dimension-dependent threshold, the outward direction is realized by
biased Bernoulli initializations. In dimension two, throughout the
coexistence regime, the resulting two-way forcing interval contains
$[-m_\beta,m_\beta]$ and extends beyond it. The complete statement is
Theorem~\ref{thm:global-forcing}.

For nonzero rational levels, the witnesses are periodic slabs.
Period-subgroup ergodicity yields a deterministic magnetization curve.
A graphical-clock estimate then converts a short-time drift bound into
pathwise exit on an interval of times. For an arbitrary nonzero interior level,
a nearby rational slab is evolved until its magnetization curve reaches
the prescribed value. The law at the hitting time is concentrated on
that level. This gives full-measure families of pathwise forcing
witnesses and a consistent semigroup family along the intermediate
levels crossed by the slab trajectory.

The centered level requires a separate construction. We use a balanced
periodic configuration with nonzero initial drift. Its cylindrical
extensions give the corresponding construction in every dimension,
and global spin-flip reverses the sign of the drift. Evolving this
configuration gives two-way forcing on an interval around zero for
every $\beta>0$.

We next consider biased Bernoulli product measures. Their initial
magnetization drift is an explicit odd polynomial in the bias. There is
a dimension-dependent inverse-temperature threshold above which this
drift is directed away from zero on a nontrivial interval. In dimension
two the drift is cubic, and its positive zero
$r_{\beta,2}^{\mathrm{Ber}}$ satisfies
\[
    r_{\beta,2}^{\mathrm{Ber}}>m_\beta
\]
throughout the coexistence regime. Hence the resulting two-way forcing
interval contains $[-m_\beta,m_\beta]$ and extends beyond it.

In dimension two, $\beta^{\mathrm{Ber}}_2<\beta_c$. For
$\beta^{\mathrm{Ber}}_2<\beta\leq\beta_c$ and
$0<r<r^{\mathrm{Ber}}_{\beta,2}$, the Bernoulli-quench magnetization
has positive initial derivative, whereas attractiveness and uniqueness imply convergence
to zero. Its magnetization curve is therefore nonmonotone.

For $\beta>\beta_c$, Proposition~\ref{prop:quench-threshold-vs-spontaneous}
gives $m_\beta<r^{\mathrm{Ber}}_{\beta,2}$. Hence, for
\[
    m_\beta<r<r^{\mathrm{Ber}}_{\beta,2},
\]
the Bernoulli-quench magnetization initially increases although
$r>m_\beta$. By attractiveness, it is bounded above by the
magnetization of the process started from the all-plus configuration,
which converges to $m_\beta$ \cite{Liggett85, Martinelli99}. The magnetization
curve is again nonmonotone. For comparison, the corresponding
long-time phase-ordering result is given in \cite{GheissariSly26}.

Consequently, nonmonotone Bernoulli-quench magnetization occurs for
every $\beta>\beta^{\mathrm{Ber}}_2$. In the coexistence regime, the
initial magnetization may exceed $m_\beta$ while retaining positive
initial drift.

Configurations in the same exact magnetization level can therefore
have opposite short-time macroscopic drifts, depending on the
frequencies of their local environments.

\section{Exact levels and pathwise forcing}
\label{sec:pathwise-forcing}

Fix an integer \(d\geq1\). We work on
\[
    \Omega=\{-1,+1\}^{\mathbb Z^d}
\]
with its product topology and Borel $\sigma$-field. For \(n\geq1\), let
\[
    \Lambda_n=[-n,n]^d\cap\mathbb Z^d,
    \qquad
    M_n(\sigma)
    =
    \frac1{|\Lambda_n|}
    \sum_{x\in\Lambda_n}\sigma(x).
\]
For \(r\in[-1,1]\), define the exact magnetization level
\[
    C_r
    =
    \left\{
        \sigma\in\Omega:
        M_n(\sigma)\longrightarrow r
    \right\}.
\]
The set \(C_r\) is Borel. It fixes the limiting density of plus and
minus spins, but does not prescribe their local arrangement.

We consider the nearest-neighbour heat-bath Glauber dynamics at inverse
temperature \(\beta>0\), with zero external field. Its generator acts
on bounded local functions according to
\[
    Lf(\sigma)
    =
    \sum_{x\in\mathbb Z^d}
    c_x(\sigma)
    \bigl[f(\sigma^x)-f(\sigma)\bigr],
\]
where \(\sigma^x\) denotes the configuration obtained from \(\sigma\)
by flipping the spin at \(x\), namely,
\[
    \sigma^x(y)
    =
    \begin{cases}
        -\sigma(x), & y=x,\\
        \sigma(y),  & y\neq x,
    \end{cases}
\]
and
\[
    c_x(\sigma)
    =
    \frac12
    \left[
        1-\sigma(x)\tanh\bigl(\beta h_x(\sigma)\bigr)
    \right],
    \qquad
    h_x(\sigma)
    =
    \sum_{y\sim x}\sigma(y).
\]
We denote the associated Markov semigroup by
\((P_t)_{t\geq0}\), and by \( \mathbb P_\xi\) the law of the process started
from \(\xi\in\Omega\).

For \(x\in\mathbb Z^d\), let \(f_x(\sigma)=\sigma(x)\). Then
\[
    Lf_x(\sigma)
    =
    \tanh\bigl(\beta h_x(\sigma)\bigr)-\sigma(x).
\]
Hence the short-time empirical drift depends on the spatial distribution
of the local fields $h_x$, not only on the initial magnetization.

\begin{definition}[Pathwise forcing]
\label{def:pathwise-forcing}
Let \(r\in[-1,1]\) and \(\xi\in C_r\).

We call \(\xi\) a \emph{downward pathwise forcing witness} for \(C_r\)
if there exists \(\varepsilon>0\) such that
\[
    \mathbb P_\xi\left(
        \limsup_{n\to\infty}M_n(\sigma_t)<r
        \text{ for every }0<t<\varepsilon
    \right)=1.
\]

We call \(\xi\) an \emph{upward pathwise forcing witness} for \(C_r\)
if there exists \(\varepsilon>0\) such that
\[
    \mathbb P_\xi\left(
        \liminf_{n\to\infty}M_n(\sigma_t)>r
        \text{ for every }0<t<\varepsilon
    \right)=1.
\]

A configuration satisfying either condition is called a
\emph{pathwise forcing witness}.
\end{definition}

In either case, the process leaves the initial exact level
simultaneously at every sufficiently small positive time. In
particular, if \(\xi\) is a pathwise forcing witness for \(C_r\), then
\[
    P_t(\xi,C_r)=0
\]
for every sufficiently small \(t>0\). The converse need not hold,
since fixed-time statements do not in general imply simultaneous
pathwise exit over an uncountable time interval.

For \(r>0\), forcing toward zero means downward forcing; for \(r<0\),
it means upward forcing. At the centered level \(C_0\), the two
directions are distinct and will both occur.

\begin{remark}[Measurability]
The path events in Definition~\ref{def:pathwise-forcing} are universally
measurable. Indeed, the maps
\[
 \sigma\longmapsto\limsup_{n\to+\infty} M_n(\sigma),
 \qquad
 \sigma\longmapsto\liminf_{n\to+\infty} M_n(\sigma)
\]
are Borel, as is the joint evaluation map
$(t,\omega)\mapsto\omega(t)$ on the Polish Skorokhod space
$D([0,\infty),\Omega)$. The complement of each forcing event is therefore
the projection onto path space of a Borel subset of
$(0,\varepsilon)\times D([0,\infty),\Omega)$, hence is analytic.
The probabilities in Definition~\ref{def:pathwise-forcing} are thus
well-defined under the completed path laws $\mathbb P_\xi$.
\end{remark}

\subsection{A pathwise thickening principle}
\label{subsec:pathwise-thickening}

The next lemma upgrades exact magnetization control on a countable
dense set of times to simultaneous control throughout a punctured
time interval. Its proof uses only the rate-one Poisson clocks in the
standard graphical construction.

\begin{lemma}[Pathwise thickening from dense exact-time control]
\label{lem:pathwise-thickening}
Let \(r\in[-1,1]\), let \(\xi\in C_r\), and let
\(\varepsilon>0\). Assume that there exist a constant \(c>0\), a
countable dense set
\[
    D\subset(0,\varepsilon),
\]
and a function \(m:D\to[-1,1]\) such that, for every \(q\in D\),
\[
    \mathbb P_\xi\left(
        M_n(\sigma_q)\longrightarrow m(q)
    \right)=1
\]
and
\[
    m(q)\le r-cq.
\]
Then \(\xi\) is a downward pathwise forcing witness for \(C_r\).
More precisely, there exists \(\varepsilon'>0\) such that
\[
    \mathbb P_\xi\left(
        \limsup_{n\to\infty}M_n(\sigma_t)<r
        \text{ for every }0<t<\varepsilon'
    \right)=1.
\]

Symmetrically, if
\[
    m(q)\ge r+cq,
    \qquad q\in D,
\]
then \(\xi\) is an upward pathwise forcing witness.
\end{lemma}

\begin{proof}
We prove the downward statement. Choose
\[
0<\alpha<\min\{1,c/4\}.
\]

We first record the elementary clock estimate. For every deterministic
\(u\geq0\) and \(h>0\),
\[
    \sup_{u\le t\le u+h}
    \limsup_{n\to\infty}
    \left|
        M_n(\sigma_t)-M_n(\sigma_u)
    \right|
    \le
    2(1-e^{-h})
    \le 2h
    \qquad \mathbb P_\xi\text{-a.s.}
\]
Indeed, if the spins at a site differ between times \(u\) and
\(t\le u+h\), then the clock at that site must have rung at least once
during \([u,u+h]\). Hence, for every such \(t\),
\[
    \left|
        M_n(\sigma_t)-M_n(\sigma_u)
    \right|
    \le
    \frac{2}{|\Lambda_n|}
    \sum_{x\in\Lambda_n}
    \mathbf 1_{\{N_x(u,u+h]\ge1\}},
\]
where \(N_x(u,u+h]\) is the number of rings of the clock at \(x\)
during \((u,u+h]\). The indicators on the right-hand side are
independent Bernoulli variables with mean \(1-e^{-h}\), so their
spatial average converges almost surely to \(1-e^{-h}\).

For each \(q\in D\), let \(E_q\) be an event of probability one on
which
\[
    M_n(\sigma_q)\longrightarrow m(q)\le r-cq
    \qquad\text{as }n\to\infty
\]
and the corresponding clock estimate over the interval
\[
    [(1-\alpha)q,(1+\alpha)q]
\]
hold simultaneously. If
\[
    |t-q|\le\alpha q,
\]
then a discrepancy between \(\sigma_t(x)\) and \(\sigma_q(x)\)
requires at least one clock ring at \(x\) during
\[
    [(1-\alpha)q,(1+\alpha)q],
\]
whose length is \(2\alpha q\). Therefore, on \(E_q\),
simultaneously for every such \(t\),
\[
\begin{aligned}
    \limsup_{n\to\infty}M_n(\sigma_t)
    \le
    m(q)+2(1-e^{-2\alpha q}) \le
    r-cq+4\alpha q
    <r.
\end{aligned}
\]
Consequently,
\[
    \limsup_{n\to\infty}M_n(\sigma_t)<r
    \qquad
    \text{for every }
    t\in I_q:=((1-\alpha)q,(1+\alpha)q).
\]

Since \(D\) is countable,
\[
    E:=\bigcap_{q\in D}E_q
\]
satisfies
\[
    \mathbb P_\xi(E)=1.
\]
Since \(D\) is dense in \((0,\varepsilon)\), the intervals
\[
    I_q=((1-\alpha)q,(1+\alpha)q),\qquad q\in D,
\]
cover \((0,(1-\alpha)\varepsilon)\). Indeed, if \(0<t<(1-\alpha)\varepsilon\), then $
\left(\frac{t}{1+\alpha},\frac{t}{1-\alpha}\right)$
is a nonempty open subinterval of \((0,\varepsilon)\), and hence contains
some \(q\in D\). For this \(q\), one has \(t\in I_q\).

Hence, on the single event \(E\),
\[
    \limsup_{n\to\infty}M_n(\sigma_t)<r
    \qquad
    \text{for every }0<t<(1-\alpha)\varepsilon.
\]
The conclusion follows with
\[
    \varepsilon'=(1-\alpha)\varepsilon.
\]

The upward statement is obtained by reversing the inequalities and
using \(\liminf\) in place of \(\limsup\).
\end{proof}

\begin{corollary}[Forcing under a drifting law]\label{forcing}
Let \(\mu\) be a probability measure such that \(\mu(C_r)=1\). Assume
that there exist \(c,\varepsilon>0\) and \(m:(0,\varepsilon)\to[-1,1]\)
such that
\[
\mu P_t(C_{m(t)})=1,
\qquad
m(t)\leq r-ct,
\qquad 0<t<\varepsilon.
\]
Then there exist a Borel set \(B\subset C_r\), with \(\mu(B)=1\), and
\(\varepsilon'>0\) such that every \(\xi\in B\) is a downward pathwise
forcing witness on \((0,\varepsilon')\).

The symmetric statement holds if \(m(t)\geq r+ct\).

\end{corollary}

\begin{proof}
Let $D\subset(0,\varepsilon)$ be countable and dense. For each $q\in D$,
\[
1=\mu P_q(C_{m(q)})
 =\int P_q(\xi,C_{m(q)})\,\mu(d\xi).
\]
Hence the Borel set
\[
B_q:=\{\xi\in C_r:\ P_q(\xi,C_{m(q)})=1\}
\]
has $\mu(B_q)=1$. Therefore
\[
B:=\bigcap_{q\in D}B_q
\]
satisfies $\mu(B)=1$.

For every $\xi\in B$ and every $q\in D$,
\[
\mathbb P_\xi\bigl(M_n(\sigma_q)\to m(q)\bigr)=1.
\]
Since $m(q)\le r-cq$, Lemma~\ref{lem:pathwise-thickening} applies to each $\xi\in B$.

\end{proof}

\subsection{The one-dimensional case}

We now dispose of the case $d=1$. Set
\[
        a_\beta:=\tanh(2\beta),
        \qquad
        \lambda_\beta:=1-a_\beta>0.
\]
The two neighbours of $x$ can have only the following configurations:
\[
\begin{array}{c|c|c}
(\sigma(x-1),\sigma(x+1))
& h_x(\sigma)
& \tanh(\beta h_x(\sigma))
\\ \hline
(+1,+1) & 2  & a_\beta\\
(+1,-1) & 0  & 0\\
(-1,+1) & 0  & 0\\
(-1,-1) & -2 & -a_\beta .
\end{array}
\]
Hence
\[
        Lf_x
        =
        \frac{a_\beta}{2}(f_{x-1}+f_{x+1})-f_x.
\]
Averaging the preceding identity over $\Lambda_n$ gives
\[
\begin{aligned}
LM_n
&=
\frac{a_\beta}{2(2n+1)}
\sum_{x=-n}^n
\bigl(f_{x-1}+f_{x+1}\bigr)
-M_n \\
&=
(a_\beta-1)M_n
+
\frac{a_\beta}{2(2n+1)}
\bigl(
f_{-n-1}+f_{n+1}-f_{-n}-f_n
\bigr).
\end{aligned}
\]
Hence
\[
LM_n=-\lambda_\beta M_n+R_n,
\]
where
\[
R_n(\sigma)
=
\frac{a_\beta}{2(2n+1)}
\bigl[
\sigma(-n-1)+\sigma(n+1)-\sigma(-n)-\sigma(n)
\bigr],
\]
and therefore
\[
\|R_n\|_\infty\le\frac{2a_\beta}{2n+1}.
\]

Fix $\xi\in C_r$. Dynkin's formula and variation of constants give
\[
M_n(\sigma_t)
=
e^{-\lambda_\beta t}M_n(\xi)
+
\int_0^t e^{-\lambda_\beta(t-s)}R_n(\sigma_s)\,ds
+
\int_0^t e^{-\lambda_\beta(t-s)}\,dZ_n(s),
\]
where $Z_n$ is a martingale. Since each jump of $M_n$ has size at most
$2/(2n+1)$ and the total update rate in $\Lambda_n$ is at most $2n+1$,
the last term has second moment $O(n^{-1})$. Taking $n_k=k^2$,
Chebyshev's inequality and Borel--Cantelli, followed by the deterministic
comparison between $M_n$ and $M_{k^2}$ for
$k^2\le n<(k+1)^2$, yield, for every fixed $t\ge0$,
\[
        M_n(\sigma_t)
        \longrightarrow
        r e^{-\lambda_\beta t}
        \qquad P_\xi\text{-a.s.}
\]
Thus the convergence holds simultaneously at all nonnegative rational
times on a single probability-one event.

If $r>0$, then, for some $c,\varepsilon>0$,
\[
        r e^{-\lambda_\beta t}\le r-ct,
        \qquad 0<t<\varepsilon.
\]
Lemma~\ref{lem:pathwise-thickening} therefore shows that every $\xi\in C_r$ is a downward pathwise
forcing witness. By spin-flip symmetry, every $\xi\in C_r$ is an upward
pathwise forcing witness when $r<0$.

If $r=0$, the empirical magnetization converges to zero at every
rational time on a single probability-one event. Intersecting this
event with the countable family of clock-estimate events corresponding
to rational time intervals, we obtain a probability-one event on which
the following holds. For every $t>0$ and every rational
$0<h<t/2$, choose $q\in\mathbb Q_+$ with $|t-q|<h$. Then
\[
    \limsup_{n\to\infty}|M_n(\sigma_t)|
    \leq 2(1-e^{-2h}).
\]
Letting $h\downarrow0$ yields
\[
    M_n(\sigma_t)\longrightarrow0
    \qquad\text{for every }t\geq0
\]
simultaneously. Hence
\[
    \mathbb P_\xi\bigl(
        \sigma_t\in C_0\ \text{for every }t\geq0
    \bigr)=1,
    \qquad \xi\in C_0.
\]
Thus, in dimension one, $C_0$ is pathwise absorbing, while every
nonzero exact level is pathwise forced toward zero.

For the remainder of the paper, we assume $d\ge2$.

\section{Forcing toward the centered level}
\label{sec:centered-forcing}

We first construct forcing witnesses directed toward zero at every
nonzero level. Rational levels admit explicit periodic witnesses,
dynamic interpolation covers the remaining interior levels, and the
extremal levels are treated directly.

\subsection{Periodic slabs}
\label{subsec:periodic-slabs}

Fix integers \(p,q\geq2\), and set $r_{p,q}
    =
    \frac{p-q}{p+q}$.
    Define the periodic slab configuration \(\xi^{p,q}\in\Omega\) by
\[
    \xi^{p,q}(i_1,\ldots,i_d)
    =
    \begin{cases}
        +1,
        &
        i_1\bmod(p+q)\in\{0,\ldots,p-1\},
        \\[1mm]
        -1,
        &
        i_1\bmod(p+q)\in\{p,\ldots,p+q-1\}.
    \end{cases}
\]
Thus one period in the first coordinate consists of \(p\) plus
hyperplanes followed by \(q\) minus hyperplanes, while the
configuration is constant in the remaining \(d-1\) coordinate
directions. It is invariant under the finite-index subgroup
\[
    \Gamma_{p,q}
    =
    ((p+q)\mathbb Z)\times\mathbb Z^{d-1}.
\]

Since one period contains \(p\) plus hyperplanes and \(q\) minus
hyperplanes,
\[
    M_n(\xi^{p,q})
    \longrightarrow
    r_{p,q}.
\]
Hence
\[
    \xi^{p,q}\in C_{r_{p,q}}.
\]

For \(t\geq0\), define the period-averaged magnetization
\[
    m_{p,q}(t)
    :=
    \frac1{p+q}
    \sum_{j=0}^{p+q-1}
    \mathbb E_{\xi^{p,q}}
    \bigl[\sigma_t(j,0,\ldots,0)\bigr].
\]

Let \(\Gamma=\Gamma_{p,q}\) and \(\xi=\xi^{p,q}\). Realize the dynamics
by independent rate-one Poisson clocks \((N_x)_{x\in\mathbb Z^d}\) and
independent update variables attached to their jump times. Denote the
whole collection of graphical marks by \(G\), and let
\[
    \Phi_t(\xi,G)
\]
be the configuration at time \(t\).

For every \(\gamma\in\Gamma\), translation covariance of the update rule
and the identity \(\tau_\gamma\xi=\xi\) imply
\[
    \Phi_t(\xi,\theta_\gamma G)
    =
    \tau_\gamma\Phi_t(\xi,G),
\]
where \(\theta_\gamma\) denotes translation of the graphical marks.
Thus the map
\[
    G\longmapsto \Phi_t(\xi,G)
\]
is measurable and \(\Gamma\)-equivariant.

The law of \(G\) is a product measure over the sites of
\(\mathbb Z^d\). Choose a set \(R\) of representatives of
\(\mathbb Z^d/\Gamma\), and, for each \(\gamma\in\Gamma\), group the
marks into the vector
\[
    \widehat G_\gamma
    :=
    \bigl(G_{\gamma+a}\bigr)_{a\in R}.
\]
Then $(\widehat G_\gamma)_{\gamma\in\Gamma}$ is an i.i.d. vector-valued
field indexed by $\Gamma$, so the action of $\Gamma$ on the graphical
noise is ergodic. Since the map
\[
        G\longmapsto \Phi_t(\xi^{p,q},G)
\]
is $\Gamma$-equivariant, its image law
\[
        \delta_{\xi^{p,q}}P_t
\]
is $\Gamma_{p,q}$-invariant and $\Gamma_{p,q}$-ergodic.
For related factor-of-i.i.d.\ and finitary coding constructions for
Markov random fields, see H\"aggstr\"om and Steif~\cite{HaggstromSteif00}.

Since $\Gamma_{p,q}$ has finite index in $\mathbb Z^d$, the spatial
ergodic theorem \cite{Walters82} may be applied separately on the finitely many cosets of
$\Gamma_{p,q}$ to give
\[
    M_n(\sigma)
    \longrightarrow
    m_{p,q}(t)
    \qquad
    \delta_{\xi^{p,q}}P_t\text{-a.s.}
\]
Equivalently,
\begin{equation}
   \mathbb P_{\xi^{p,q}}
    \left(
        M_n(\sigma_t)\longrightarrow m_{p,q}(t)
    \right)
    =1,
    \qquad t\geq0.
    \label{eq:periodic-ergodic-limit}
\end{equation}
In particular,
\[
    m_{p,q}(0)=r_{p,q}.
\]

\begin{proposition}[Slab drift]
\label{prop:slab-drift}
For every \(p,q\geq2\),
\[
    m'_{p,q}(0)
    =
    r_{p,q}
    \bigl(\tanh(2d\beta)-1\bigr).
\]
Consequently, whenever \(r_{p,q}\neq0\), the initial slab drift is
directed toward zero.
\end{proposition}

\begin{proof}
Differentiating the period-averaged magnetization at time zero and
using the generator formula gives
\[
    m'_{p,q}(0)
    =
    \frac1{p+q}
    \sum_{j=0}^{p+q-1}
    \left[
        \tanh\!\bigl(\beta h_{(j,0,\ldots,0)}
        (\xi^{p,q})\bigr)
        -
        \xi^{p,q}(j,0,\ldots,0)
    \right].
\]
There are four types of hyperplanes, according to the spin and
whether the hyperplane lies in the interior or at the boundary of a
slab. Their respective local fields and contributions are
\[
\renewcommand{\arraystretch}{1.15}
\begin{array}{c@{\qquad}c@{\qquad}c@{\qquad}c}
    \text{type}
    &
    \text{multiplicity}
    &
    \text{local field}
    &
    \text{contribution}
    \\ \hline
    \text{interior plus}
    &
    p-2
    &
    2d
    &
    \tanh(2d\beta)-1
    \\
    \text{boundary plus}
    &
    2
    &
    2d-2
    &
    \tanh((2d-2)\beta)-1
    \\
    \text{interior minus}
    &
    q-2
    &
    -2d
    &
    1-\tanh(2d\beta)
    \\
    \text{boundary minus}
    &
    2
    &
    -(2d-2)
    &
    1-\tanh((2d-2)\beta)
\end{array}
\]
The two boundary contributions cancel exactly. Hence, averaging the
generator over one period cell gives
\[
\begin{aligned}
    m'_{p,q}(0)
    =
    \frac{
        (p-2)\bigl(\tanh(2d\beta)-1\bigr)
        +(q-2)\bigl(1-\tanh(2d\beta)\bigr)
    }{p+q}    
    =
    \frac{p-q}{p+q}
    \bigl(\tanh(2d\beta)-1\bigr),
\end{aligned}
\]
as claimed.

\end{proof}

\begin{corollary}[Pathwise slab forcing at nonzero rational levels]
\label{cor:pathwise-rational-forcing}
For every
\[
    r\in(-1,1)\cap\mathbb Q\setminus\{0\},
\]
there exist a periodic configuration
\[
    \xi_r^{\mathrm{slab}}\in C_r
\]
and a number \(\delta_r>0\) such that \(\xi_r^{\mathrm{slab}}\) is a
pathwise forcing witness directed toward zero.

More precisely, if \(r>0\), then
\[
   \mathbb P_{\xi_r^{\mathrm{slab}}}
    \left(
        \limsup_{n\to\infty}M_n(\sigma_t)<r
        \text{ for every }0<t<\delta_r
    \right)
    =1,
\]
whereas if \(r<0\), then
\[
    \mathbb P_{\xi_r^{\mathrm{slab}}}
    \left(
        \liminf_{n\to\infty}M_n(\sigma_t)>r
        \text{ for every }0<t<\delta_r
    \right)
    =1.
\]
\end{corollary}

\begin{proof} Assume first that \(1>r>0\).
If \(r=a/b\), one may take
\[
p=N(b+a),\qquad q=N(b-a),
\]
where \(N\) is chosen large enough that \(p,q\ge2\).
Set
\[
\xi_r^{\mathrm{slab}}:=\xi^{p,q}.
\]
Then \(\xi_r^{\mathrm{slab}}\in C_r\) and
\[
m_{p,q}(0)=r.
\]
By Proposition~\ref{prop:slab-drift},
\[
m'_{p,q}(0)
=
r\bigl(\tanh(2d\beta)-1\bigr)<0.
\]
Since \(L f_x\) is a bounded local function, the map
\[
t\longmapsto P_tLf_x(\xi_r^{\mathrm{slab}})
\]
is continuous. Averaging over one period shows that \(m'_{p,q}\) is
continuous at zero. Hence there exist constants \(c>0\) and
\(\varepsilon>0\) such that
\[
m'_{p,q}(t)\le -c,
\qquad 0\le t\le\varepsilon.
\]
Therefore
\[
m_{p,q}(t)\le r-ct,
\qquad 0\le t\le\varepsilon.
\]

Moreover, by~(\ref{eq:periodic-ergodic-limit}),
\[
\delta_{\xi_r^{\mathrm{slab}}}P_t
\bigl(C_{m_{p,q}(t)}\bigr)=1,
\qquad t\ge0,
\]
while
\[
\delta_{\xi_r^{\mathrm{slab}}}(C_r)=1.
\]
Corollary~\ref{forcing}, applied with
\[
\mu=\delta_{\xi_r^{\mathrm{slab}}},
\qquad
m(t)=m_{p,q}(t),
\]
therefore shows that
\(\xi_r^{\mathrm{slab}}\) is a downward pathwise forcing witness for
\(C_r\).

The case \(-1<r<0\) follows by global spin-flip symmetry, which sends
\(C_r\) to \(C_{-r}\) and reverses the direction of forcing.
\end{proof}

\subsection{Uniform persistence of the slab drift}

To transport rational slab constructions to arbitrary nonzero interior
levels, we need the sign of the slab drift to persist uniformly for a
short time when the initial magnetization stays away from zero.

\begin{lemma}[Uniform short-time control of the slab drift]
\label{lem:uniform-drift-control}
Set
\[
    C_{\beta,d}:=1+\|L^2f_0\|_\infty<\infty.
\]
Then, for every $x\in\mathbb Z^d$, every $\xi\in\Omega$, and every
$t\ge0$,
\[
    \bigl|P_tLf_x(\xi)-Lf_x(\xi)\bigr|
    \le C_{\beta,d}t.
\]
Moreover, for every $a\in(0,1)$, there exist constants
$c_{a,\beta,d}>0$ and $t_{a,\beta,d}>0$ such that, for every rational
$\rho\in[a,1)$, the magnetization curve $m_\rho$ associated with the
periodic slab configuration of level $\rho$ satisfies
\[
    m_\rho'(t)\le -c_{a,\beta,d},
    \qquad 0\le t\le t_{a,\beta,d}.
\]
\end{lemma}

\begin{proof}
The function $Lf_0$ is bounded and local, and hence so is $L^2f_0$.
Translation covariance gives
\[
    \|L^2f_x\|_\infty=\|L^2f_0\|_\infty<C_{\beta,d},
    \qquad x\in\mathbb Z^d.
\]
Dynkin's formula applied to $Lf_x$ therefore yields
\[
\begin{aligned}
    \bigl|P_tLf_x(\xi)-Lf_x(\xi)\bigr|
    =
    \left|\int_0^t P_sL^2f_x(\xi)\,ds\right|  
    \le t\|L^2f_x\|_\infty
    \le C_{\beta,d}t.
\end{aligned}
\]

Set
\[
    \kappa_{\beta,d}:=1-\tanh(2d\beta)>0.
\]
For a rational $\rho\in(0,1)$,
Proposition~\ref{prop:slab-drift} gives
\[
    m_\rho'(0)=-\rho\kappa_{\beta,d}.
\]
Averaging the preceding estimate over one period of
$\xi_\rho^{\mathrm{slab}}$ gives
\[
    |m_\rho'(t)-m_\rho'(0)|
    \le C_{\beta,d}t.
\]

Now fix $a\in(0,1)$ and set
\[
    c_{a,\beta,d}:=\frac{a\kappa_{\beta,d}}2,
    \qquad
    t_{a,\beta,d}:=\frac{a\kappa_{\beta,d}}{2C_{\beta,d}}.
\]
For every rational $\rho\in[a,1)$ and every
$0\le t\le t_{a,\beta,d}$,
\[
\begin{aligned}
    m_\rho'(t)
    \le m_\rho'(0)+C_{\beta,d}t  \le -a\kappa_{\beta,d}
          +\frac{a\kappa_{\beta,d}}2
     =-c_{a,\beta,d}.
\end{aligned}
\]
\end{proof}

\subsection{Dynamic interpolation and forcing under interpolated laws}

The semigroup identity provides exact magnetization control at each
fixed time. Restricting to a countable dense set of times and using
Lemma~\ref{lem:pathwise-thickening} then thickens this control into a
whole punctured interval. A single periodic slab trajectory therefore
generates full-measure families of pathwise forcing witnesses at every
intermediate magnetization level it crosses. In particular, the
construction applies without distinction to rational and irrational
nonzero interior levels.

\begin{theorem}[Full-measure forcing under an interpolated law] \label{thm:typical-toward-zero}
For every \(r\in(-1,1)\setminus\{0\}\), there exist a probability
measure \(\nu_r\) concentrated on \(C_r\), a Borel set
\(B_r\subset C_r\), and \(\varepsilon_r>0\) such that
\[
    \nu_r(B_r)=1,
\]
and every \(\xi\in B_r\) is a pathwise forcing witness directed
toward zero on \((0,\varepsilon_r)\).
\end{theorem}

\begin{proof}
We treat the case \(r>0\); the case \(r<0\) follows by global
spin-flip symmetry.

Apply Lemma~\ref{lem:uniform-drift-control} with \(a=r\), and denote
the resulting constants by \(c_{r,\beta,d}>0\) and
\(t_{r,\beta,d}>0\). Choose
\[
    \rho\in(r,1)\cap\mathbb Q
\]
so close to \(r\) that
\[
    0<\rho-r<c_{r,\beta,d}t_{r,\beta,d}.
\]
Let \(\xi_\rho^{\mathrm{slab}}\) be the periodic slab configuration of
level \(\rho\), and let \(m_\rho\) be its deterministic magnetization
curve.

Since
\[
    m_\rho(0)=\rho>r
\]
and
\[
    m_\rho(t_{r,\beta,d})
    \le
    \rho-c_{r,\beta,d}t_{r,\beta,d}
    <r,
\]
continuity of \(m_\rho\) yields \(s\in(0,t_{r,\beta,d})\) such that
\[
    m_\rho(s)=r.
\]
Set
\[
    \nu_r:=\delta_{\xi_\rho^{\mathrm{slab}}}P_s.
\]
Since \(\xi_\rho^{\mathrm{slab}}\) is periodic, \(\nu_r\) is invariant
and ergodic under its period subgroup. Hence the spatial ergodic
theorem gives
\[
    M_n(\sigma)\longrightarrow m_\rho(s)=r
    \qquad
    \nu_r\text{-a.s.},
\]
and therefore
\[
    \nu_r(C_r)=1.
\]

Choose \(\varepsilon_r>0\) such that
\[
    s+\varepsilon_r<t_{r,\beta,d}.
\]
Since
\[
    m_\rho'(u)\le -c_{r,\beta,d},
    \qquad 0\le u\le t_{r,\beta,d},
\]
we have, for every \(0\le t\le\varepsilon_r\),
\[
\begin{aligned}
    m_\rho(s+t)-m_\rho(s)
    =
    \int_s^{s+t}m_\rho'(u)\,du 
    \le
    -c_{r,\beta,d}t.
\end{aligned}
\]
Thus
$$
    m_\rho(s+t)\le r-c_{r,\beta,d}t,
    \qquad 0\le t\le\varepsilon_r.
$$

For every \(t\in(0,\varepsilon_r)\), the semigroup property gives
\[
    \nu_rP_t
    =
    \delta_{\xi_\rho^{\mathrm{slab}}}P_{s+t}.
\]
Again by periodicity and the spatial ergodic theorem,
$$
    \nu_rP_t\bigl(C_{m_\rho(s+t)}\bigr)=1.
$$

The preceding two relations, together with
Corollary~\ref{forcing}, applied with
\[
    \mu=\nu_r,
    \qquad
    m(t)=m_\rho(s+t),
\]
yield a Borel set \(B_r\subset C_r\), with \(\nu_r(B_r)=1\), whose
elements are downward pathwise forcing witnesses.
This proves the result for \(r>0\). The case \(r<0\) follows by global
spin-flip symmetry.
\end{proof}

Fix \(0<r<\rho<1\), where \(\rho\) is rational and chosen as in the
proof of Theorem~\ref{thm:typical-toward-zero}. Let \(m_\rho\) be the magnetization curve of
\(\xi_\rho^{\mathrm{slab}}\). For every \(u\in[r,\rho]\), let
\(s_u\in[0,t_{r,\beta,d})\) be the unique time such that
\[
    m_\rho(s_u)=u,
\]
and set
\[
    \nu_u:=\delta_{\xi_\rho^{\mathrm{slab}}}P_{s_u}.
\]

\begin{corollary}[Consistent descendant forcing laws]
\label{cor:descendant-forcing-laws}
There exist a common time $\varepsilon_{\mathrm{pw}}>0$ and, for every
$u\in[r,\rho]$, a Borel set $B_u\subset C_u$ such that
\[
    \nu_u(B_u)=1,
\]
and every $\xi\in B_u$ satisfies
\[
    \mathbb P_\xi\left(
        \limsup_{n\to\infty}M_n(\sigma_t)<u
        \text{ for every }0<t<\varepsilon_{\mathrm{pw}}
    \right)=1.
\]

Moreover, for every $r\le u<v\le\rho$,
\[
    \nu_vP_{s_u-s_v}=\nu_u,
\]
and
\[
    P_{s_u-s_v}(\eta,B_u)=1
    \qquad\text{for $\nu_v$-almost every $\eta\in C_v$.}
\]
Thus the descendant law
\[
    \delta_\eta P_{s_u-s_v}
\]
is concentrated on $B_u$ for $\nu_v$-almost every $\eta$.
\end{corollary}

\begin{proof}
Set
\[
    \varepsilon_0:=t_{r,\beta,d}-s_r>0.
\]
Since $s_u\le s_r$ for every $u\in[r,\rho]$, the uniform drift bound
gives
\[
    m_\rho(s_u+t)\le u-c_{r,\beta,d}t,
    \qquad
    0\le t\le\varepsilon_0.
\]

Moreover,
\[
    \nu_uP_t
    =
    \delta_{\xi_\rho^{\mathrm{slab}}}P_{s_u+t},
\]
and periodicity implies
\[
    \nu_uP_t\bigl(C_{m_\rho(s_u+t)}\bigr)=1.
\]
The proof of Corollary~\ref{forcing}, applied with the same constants
$c_{r,\beta,d}$ and $\varepsilon_0$ for every $u$, yields a common
$\varepsilon_{\mathrm{pw}}>0$ and the sets $B_u$.

For $r\le u<v\le\rho$, the semigroup property gives
\[
    \nu_vP_{s_u-s_v}=\nu_u.
\]
Since $\nu_u(B_u)=1$,
\[
    \int P_{s_u-s_v}(\eta,B_u)\,\nu_v(d\eta)=1,
\]
and hence
\[
    P_{s_u-s_v}(\eta,B_u)=1
\]
for $\nu_v$-almost every $\eta$.
\end{proof}

For a fixed target level $u$, the descendant laws
$\delta_\eta P_{s_u-s_v}$, $v\in(u,\rho]$, are concentrated on
$B_u$ for $\nu_v$-almost every $\eta$. They are consistent under
composition: if $r\leq u<v<w\leq\rho$, then
\[
    \delta_\eta P_{s_u-s_w}
    =
    \bigl(\delta_\eta P_{s_v-s_w}\bigr)P_{s_u-s_v}.
\]

In particular, the rational periodic slab is only used to start the
construction. Every intermediate magnetization level, rational or
irrational, carries a full-measure family of pathwise forcing witnesses.
Moreover, the construction can be restarted from almost every
configuration drawn at an intermediate level: evolution to any lower
level produces a forcing witness there almost surely. The procedure may
be iterated through any finite chain of intermediate levels.

\begin{remark}[Explicit aperiodic candidates at irrational levels]
For irrational \(r\in(-1,1)\), a Sturmian slab configuration with
plus-spin density
\[
    \alpha=\frac{1+r}{2}
\]
belongs to \(C_r\); see~\cite{lothaire02}. Proving that it is a forcing
witness would additionally require determining the sign of its averaged
local-environment drift. Theorem~\ref{thm:typical-toward-zero} instead provides a full-measure
family of forcing witnesses at level \(r\) without selecting an
explicit representative.
\end{remark}

\subsection{The extremal levels} \label{extreme}

It remains to treat the levels \(r=\pm1\). Let \(\xi^+\) denote the
all-plus configuration. Then
\[
    \xi^+\in C_1.
\]
For \(t\geq0\), let
\[
    m^+(t):=\mathbb E_{\xi^+}[\sigma_t(0)].
\]
Since the initial configuration is translation invariant, the law
\(\delta_{\xi^+}P_t\) is translation invariant and ergodic. Hence
\[
    \mathbb P_{\xi^+}\bigl(M_n(\sigma_t)\longrightarrow m^+(t)\bigr)=1,
    \qquad t\geq0.
\]
Moreover,
\[
    m^+(0)=1,
    \qquad
    (m^+)'(0)=\tanh(2d\beta)-1<0.
\]
By continuity of \((m^+)'\) at zero, there exist \(c,\varepsilon>0\)
such that
\[
    m^+(t)\leq1-ct,
    \qquad 0\leq t\leq\varepsilon.
\]
Let \(D\subset(0,\varepsilon)\) be countable and dense. For every
\(q\in D\),
\[
\mathbb    P_{\xi^+}\bigl(M_n(\sigma_q)\longrightarrow m^+(q)\bigr)=1,
    \qquad
    m^+(q)\leq1-cq.
\]
Lemma~\ref{lem:pathwise-thickening} shows that \(\xi^+\) is a downward
pathwise forcing witness for \(C_1\). The case \(r=-1\) follows by
global spin-flip symmetry.

\section{Biased quenches and opposite pathwise forcing}
\label{sec:biased-quenches}

The centered construction gives a possibly nonexplicit two-way forcing
interval around zero. We now consider outward forcing under biased
Bernoulli product measures.
Its drift is explicit in every dimension and, in dimension
two, can be compared with the spontaneous magnetization.

\subsection{Biased Bernoulli quenches}

For \(r\in[-1,1]\), let \(\pi_r\) denote the Bernoulli product measure
on \(\Omega\) defined by
\[
    \pi_r(\sigma(x)=+1)
    =
    \frac{1+r}{2},
    \qquad
    \pi_r(\sigma(x)=-1)
    =
    \frac{1-r}{2}.
\]
The measure \(\pi_r\) is invariant and ergodic under the full
\(\mathbb Z^d\)-action, and
\[
    \int_\Omega \sigma(0)\,d\pi_r=r.
\]
Hence, by the pointwise ergodic theorem,
\[
    \pi_r(C_r)=1.
\]

We denote by
\[
    m_r^{\mathrm{Ber}}(t)
    :=
    \int_\Omega P_tf_0\,d\pi_r
\]
the deterministic magnetization curve of the biased product quench.
Thus
\[
    m_r^{\mathrm{Ber}}(0)=r,
    \qquad
    \bigl(m_r^{\mathrm{Ber}}\bigr)'(0)
    =
    D_{\beta,d}^{\mathrm{Ber}}(r),
\]
where
\[
    D_{\beta,d}^{\mathrm{Ber}}(r)
    :=
    \int_\Omega Lf_0\,d\pi_r
\]
is its initial magnetization drift.

Let \(Y_1,\ldots,Y_{2d}\) be independent random variables with
\[
    \mathbb P(Y_j=+1)=\frac{1+r}{2},
    \qquad
    \mathbb P(Y_j=-1)=\frac{1-r}{2},
\]
and set
\[
    S_{2d}=Y_1+\cdots+Y_{2d}.
\]
Under \(\pi_r\), the local field at the origin has the same law as
\(S_{2d}\).

\begin{proposition}[Drift of the biased quench]
\label{prop:quench-drift-general}
For every \(d\geq2\), \(\beta>0\), and \(r\in[-1,1]\),
\[
    D_{\beta,d}^{\mathrm{Ber}}(r)
    =
    \mathbb E_{\mathrm{Ber}(r)}
    \left[
        \tanh(\beta S_{2d})
    \right]
    -
    r.
\]
Equivalently,
\[
D_{\beta,d}^{\mathrm{Ber}}(r)
=
\sum_{k=0}^{2d}
\binom{2d}{k}
\left(\frac{1+r}{2}\right)^k
\left(\frac{1-r}{2}\right)^{2d-k}
\tanh\!\bigl(2\beta(k-d)\bigr)
-r.
\]
The function \(r\mapsto D_{\beta,d}^{\mathrm{Ber}}(r)\) is an odd
polynomial of degree at most \(2d-1\).
\end{proposition}

\begin{proof}
By the heat-bath generator formula,
\[
    D_{\beta,d}^{\mathrm{Ber}}(r)
    =
    \int_\Omega
    \left[
        \tanh\!\bigl(\beta h_0(\sigma)\bigr)-\sigma(0)
    \right]
    \pi_r(d\sigma).
\]
Under \(\pi_r\), the spin at the origin has mean \(r\), while the local
field \(h_0\) has the same law as \(S_{2d}\). Therefore
\[
    D_{\beta,d}^{\mathrm{Ber}}(r)
    =
    \mathbb E_{\mathrm{Ber}(r)}[\tanh(\beta S_{2d})]-r.
\]

If exactly \(k\) among the \(2d\) neighbouring spins are positive,
then
\[
    S_{2d}=2(k-d),
\]
and this event has probability
\[
    \binom{2d}{k}
    \left(\frac{1+r}{2}\right)^k
    \left(\frac{1-r}{2}\right)^{2d-k}.
\]
This gives the explicit binomial formula.

Finally, global spin-flip sends \(\pi_r\) to \(\pi_{-r}\) and changes
the sign of the local field. Hence
\[
    D_{\beta,d}^{\mathrm{Ber}}(-r)
    =
    -D_{\beta,d}^{\mathrm{Ber}}(r).
\]
The binomial expression shows that the drift is a polynomial of degree
at most \(2d\). Since it is odd, its degree is at most \(2d-1\).
\end{proof}

\subsection{Onset of outward Bernoulli-quench drift}

Since the quench drift is odd,
\[
    D_{\beta,d}^{\mathrm{Ber}}(0)=0.
\]
Its behaviour near the centered level is determined by the following
linear-response formula.

\begin{proposition}[Linear response at the centered quench]
\label{prop:quench-linear-response}
For every \(d\geq2\),
\[
    \left(D_{\beta,d}^{\mathrm{Ber}}\right)'(0)
    =
    \mathbb E_{\mathrm{Ber}(0)}
    \left[
        S_{2d}\tanh(\beta S_{2d})
    \right]
    -1,
\]
where, under \(\mathbb E_{\mathrm{Ber}(0)}\), the variables
\(Y_1,\ldots,Y_{2d}\) are independent symmetric spins.

Moreover, the function
\[
    \beta\longmapsto
    \mathbb E_{\mathrm{Ber}(0)}
    \left[
        S_{2d}\tanh(\beta S_{2d})
    \right]
\]
is strictly increasing from \(0\) to
\(\mathbb E_{\mathrm{Ber}(0)}|S_{2d}|\). Since
\[
    \mathbb E_{\mathrm{Ber}(0)}|S_{2d}|>1,
    \qquad d\geq2,
\]
there exists a unique number \(\beta_d^{\mathrm{Ber}}>0\) such that
\[
    \mathbb E_{\mathrm{Ber}(0)}
    \left[
        S_{2d}\tanh(\beta_d^{\mathrm{Ber}}S_{2d})
    \right]
    =1.
\]
Consequently,
\[
    \left(D_{\beta,d}^{\mathrm{Ber}}\right)'(0)>0
    \quad\Longleftrightarrow\quad
    \beta>\beta_d^{\mathrm{Ber}}.
\]
\end{proposition}

\begin{proof} For every function $F$ of the neighbouring spins,
\[
    \mathbb E_{\mathrm{Ber}(r)}[F(Y_1,\ldots,Y_{2d})]
    =
    2^{-2d}
    \sum_{y\in\{-1,1\}^{2d}}
    F(y)\prod_{j=1}^{2d}(1+ry_j).
\]
Differentiating at $r=0$ gives
\[
    \left.
    \frac{d}{dr}
    \mathbb E_{\mathrm{Ber}(r)}[F(Y_1,\ldots,Y_{2d})]
    \right|_{r=0}
    =
    2^{-2d}
    \sum_{y\in\{-1,1\}^{2d}}
    F(y)\sum_{j=1}^{2d}y_j
    =
    \mathbb E_{\mathrm{Ber}(0)}\!\left[S_{2d}F(Y_1,\ldots,Y_{2d})\right].
\]
Taking $F=\tanh(\beta S_{2d})$ therefore gives
\[
    \left(D_{\beta,d}^{\mathrm{Ber}}\right)'(0)
    =
    \mathbb E_{\mathrm{Ber}(0)}
    \left[
        S_{2d}\tanh(\beta S_{2d})
    \right]
    -1.
\]

Moreover,
\[
    \frac{d}{d\beta}
    \mathbb E_{\mathrm{Ber}(0)}
    \left[
        S_{2d}\tanh(\beta S_{2d})
    \right]
    =
    \mathbb E_{\mathrm{Ber}(0)}
    \left[
        \frac{S_{2d}^2}{\cosh^2(\beta S_{2d})}
    \right]
    >0.
\]
The function is therefore strictly increasing. It vanishes at
\(\beta=0\), and dominated convergence gives
\[
    S_{2d}\tanh(\beta S_{2d})
    \longrightarrow
    |S_{2d}|
\]
as \(\beta\to\infty\).

Finally,
\[
    \mathbb E_{\mathrm{Ber}(0)}|S_{2d}|
    =
    \frac{2d}{2^{2d}}\binom{2d}{d}.
\]
This quantity equals \(3/2\) when \(d=2\) and is strictly increasing
in \(d\), since the ratio of its values at \(d+1\) and \(d\) is
\[
    \frac{2d+1}{2d}>1.
\]
Hence
\[
    \mathbb E_{\mathrm{Ber}(0)}|S_{2d}|>1,
    \qquad d\geq2,
\]
which proves the existence and uniqueness of
\(\beta_d^{\mathrm{Ber}}\).
\end{proof}

Assume now that
\[
    \beta>\beta_d^{\mathrm{Ber}}.
\]
By Proposition~\ref{prop:quench-linear-response},
\[
    D_{\beta,d}^{\mathrm{Ber}}(0)=0,
    \qquad
    \left(D_{\beta,d}^{\mathrm{Ber}}\right)'(0)>0.
\]
The quench drift is therefore positive on some nontrivial interval to
the right of the origin.

We define the right endpoint of the maximal positive-drift interval
adjacent to the origin by
\[
    r_{\beta,d}^{\mathrm{Ber}}
    :=
    \sup
    \left\{
        \rho\in(0,1]:
        D_{\beta,d}^{\mathrm{Ber}}(r)>0
        \text{ for every }r\in(0,\rho)
    \right\}.
\]
The preceding linear-response result implies that
\[
    r_{\beta,d}^{\mathrm{Ber}}>0.
\]
Moreover,
\[
    D_{\beta,d}^{\mathrm{Ber}}(r)>0,
    \qquad
    0<r<r_{\beta,d}^{\mathrm{Ber}}.
\]

Since
\[
    D_{\beta,d}^{\mathrm{Ber}}(1)
    =
    \tanh(2d\beta)-1<0,
\]
one also has
\[
    r_{\beta,d}^{\mathrm{Ber}}<1.
\]
By continuity,
\[
    D_{\beta,d}^{\mathrm{Ber}}
    \left(r_{\beta,d}^{\mathrm{Ber}}\right)
    =0.
\]
In arbitrary dimension, we do not need to determine whether this is
the only positive zero or whether the drift may become positive again
at larger magnetizations.

\begin{theorem}[Two-way pathwise forcing with a Bernoulli-quench branch] \label{thm:two-way-general-d}
Let \(d\ge2\) and \(\beta>\beta_d^{\mathrm{Ber}}\). For every
\(0<|r|<r_{\beta,d}^{\mathrm{Ber}}\), the level \(C_r\) contains pathwise
forcing witnesses in both directions. 
For \(r>0\), downward forcing holds on a full-measure subset of \(C_r\)
under an interpolated law \(\nu_r\), whereas upward forcing holds on a
full-measure subset under the Bernoulli product law \(\pi_r\).
The statement for \(r<0\) follows by spin-flip symmetry.
\end{theorem}

\begin{proof}
Let \(0<r<r_{\beta,d}^{\mathrm{Ber}}\). By
Theorem~\ref{thm:typical-toward-zero}, there exist a probability
measure \(\nu_r\) concentrated on \(C_r\) and a Borel set
\(B_r^\downarrow\subset C_r\) with
\[
    \nu_r(B_r^\downarrow)=1,
\]
such that every \(\xi\in B_r^\downarrow\) is a downward pathwise
forcing witness.

Since
\[
    D_{\beta,d}^{\mathrm{Ber}}(r)>0,
\]
continuity of \((m_r^{\mathrm{Ber}})'\) yields \(c,\varepsilon>0\) such that
\[
    m_r^{\mathrm{Ber}}(t)\ge r+ct,
    \qquad 0\le t\le\varepsilon.
\]

Since \(\pi_rP_t\) is an equivariant factor of the product of the
Bernoulli initial field and the independent graphical noise, it is
translation invariant and ergodic. Hence
\[
\pi_rP_t\bigl(C_{m_r^{\mathrm{Ber}}(t)}\bigr)=1,
\qquad t\ge0.
\] 
The upward part of Corollary~\ref{forcing} therefore yields a Borel set
\(B_r^\uparrow\subset C_r\), with
\[
    \pi_r(B_r^\uparrow)=1,
\]
whose elements are upward pathwise forcing witnesses.

Thus \(C_r\) contains pathwise forcing witnesses in both directions,
with downward forcing holding on a \(\nu_r\)-full-measure set and
upward forcing holding on a \(\pi_r\)-full-measure set.

\end{proof}

\subsection{The exact drift-reversal level in dimension two}

In arbitrary dimension, \(r_{\beta,d}^{\mathrm{Ber}}\) is defined as the
endpoint of the first positive-drift interval issuing from the origin.
In dimension two, the quench drift is cubic and its unique positive
zero can be computed explicitly and compared with the spontaneous
magnetization.

We now specialize to \(d=2\) and set
\[
    a_\beta=\tanh(4\beta),
    \qquad
    b_\beta=\tanh(2\beta).
\]

\begin{proposition}[Explicit two-dimensional quench drift]
\label{prop:quench-drift-d2}
For $d=2$,
\[
    D_{\beta,2}^{\mathrm{Ber}}(r)
    =
    r
    \left[
        \frac{a_\beta}{2}
        +
        b_\beta
        -
        1
        +
        \left(
            \frac{a_\beta}{2}-b_\beta
        \right)r^2
    \right].
\]
Moreover,
\[
    b_\beta-\frac{a_\beta}{2}>0.
\]
Whenever
\[
    \frac{a_\beta}{2}+b_\beta-1>0,
\]
define
\[
    \left(r_{\beta,2}^{\mathrm{Ber}}\right)^2
    =
    \frac{
        \frac12a_\beta+b_\beta-1
    }{
        b_\beta-\frac12a_\beta
    }.
\]
Then
\[
    D_{\beta,2}^{\mathrm{Ber}}(r)>0
    \qquad
    \text{for }0<r<r_{\beta,2}^{\mathrm{Ber}},
\]
and
\[
    D_{\beta,2}^{\mathrm{Ber}}(r)<0
    \qquad
    \text{for }r_{\beta,2}^{\mathrm{Ber}}<r\leq1.
\]
\end{proposition}

\begin{proof}
Under $\pi_r$, the local field takes values
\[
    -4,-2,0,2,4.
\]
Using the binomial probabilities of these five values and the oddness
of the hyperbolic tangent gives
\[
    \mathbb E_{\mathrm{Ber}(r)}[\tanh(\beta S_4)]
    =
    \frac{a_\beta}{2}r(1+r^2)
    +
    b_\beta r(1-r^2).
\]
Subtracting $r$ yields
\[
    D_{\beta,2}^{\mathrm{Ber}}(r)
    =
    r
    \left[
        \frac{a_\beta}{2}
        +
        b_\beta
        -
        1
        +
        \left(
            \frac{a_\beta}{2}-b_\beta
        \right)r^2
    \right].
\]

If
\[
    u=\tanh(2\beta),
\]
then
\[
    a_\beta
    =
    \tanh(4\beta)
    =
    \frac{2u}{1+u^2}.
\]
Therefore
\[
    b_\beta-\frac{a_\beta}{2}
    =
    u-\frac{u}{1+u^2}
    =
    \frac{u^3}{1+u^2}
    >0.
\]

Thus the bracket is strictly decreasing as a function of $r^2$.
For $\beta>\beta_2^{\mathrm{Ber}}$, its value at $r=0$ is positive,
and its unique positive zero satisfies
\[
    \left(r_{\beta,2}^{\mathrm{Ber}}\right)^2
    =
    \frac{\frac12 a_\beta+b_\beta-1}
         {b_\beta-\frac12 a_\beta}.
\]
Since $a_\beta<1$, this zero belongs to $(0,1)$. The sign statement
follows.
\end{proof}

\begin{remark}[Bernoulli outward forcing in the uniqueness regime]
Let
\[
\beta_c=\frac12\log(1+\sqrt2)
\]
be the critical inverse temperature of the two-dimensional nearest-neighbour
Ising model. Since
\[
\beta\longmapsto
\frac12\tanh(4\beta)+\tanh(2\beta)-1
\]
is strictly increasing and positive at $\beta_c$, one has
\[
\beta^{\mathrm{Ber}}_2<\beta_c.
\]
More explicitly, if $u_\star\in(0,1)$ denotes the unique root of
\[
u^3-u^2+2u-1=0,
\]
then
\[
\beta^{\mathrm{Ber}}_2
=
\frac12\operatorname{arctanh}(u_\star)
\simeq 0.323643,
\qquad
\beta_c
=
\frac12\log(1+\sqrt2)
\simeq 0.440687.
\]

Hence, for every
\[
\beta^{\mathrm{Ber}}_2<\beta\leq\beta_c,
\]
there is a nontrivial interval of nonzero magnetization levels on which the
outward forcing direction is realized by biased Bernoulli quenches, although
the model is still in the uniqueness regime. Thus outward Bernoulli-quench forcing occurs strictly inside the
uniqueness regime. 
The resulting nonmonotonicity of the Bernoulli-quench magnetization curve is
discussed in Remark~\ref{rmk:quench-monotonicity}.
\end{remark}

\subsection{Two-dimensional overshoot beyond \(m_\beta\)}

Assume now that
\[
    \beta>\beta_c
\]
for the two-dimensional nearest-neighbour Ising model. Its spontaneous
magnetization is
\[
    m_\beta
    =
    \left[
        1-\sinh(2\beta)^{-4}
    \right]^{1/8}.
\]

If \(u=\tanh(2\beta)\), then
\[
    \frac{a_\beta}{2}+b_\beta-1
    =
    \frac{u}{1+u^2}+u-1.
\]
This function is strictly increasing on \((0,1)\), and at
\(u=1/\sqrt2\) its value is \(5\sqrt2/6-1>0\). Hence the condition in
Proposition~\ref{prop:quench-drift-d2} holds throughout the
coexistence regime.

\begin{proposition}[The drift-reversal level exceeds the coexistence interval]
\label{prop:quench-threshold-vs-spontaneous}
For every finite
\[
    \beta>\beta_c,
\]
one has
\[
    m_\beta<r_{\beta,2}^{\mathrm{Ber}}<1.
\]
\end{proposition}

\begin{proof}
Set
\[
    u=\tanh(2\beta).
\]
In the phase-coexistence regime,
\[
    \frac1{\sqrt2}<u<1.
\]
Using
\[
    \tanh(4\beta)=\frac{2u}{1+u^2},
\]
the drift-reversal level becomes
\[
    \bigl(r_{\beta,2}^{\mathrm{Ber}}\bigr)^2
    =
    \frac{u^3-u^2+2u-1}{u^3}
    =
    1-\frac{(1-u)^2}{u^3}.
\]
In particular,
\[
    r_{\beta,2}^{\mathrm{Ber}}<1.
\]

For the two-dimensional nearest-neighbor Ising model, the spontaneous
magnetization is given by Yang's formula~\cite{Yang52}
\[
m_\beta
=
\left(1-\sinh(2\beta)^{-4}\right)^{1/8},
\qquad \beta>\beta_c.
\]
Therefore
\[
\begin{aligned}
    \bigl(r_{\beta,2}^{\mathrm{Ber}}\bigr)^8-m_\beta^8
    =
    \left(
        \frac{u^3-u^2+2u-1}{u^3}
    \right)^4
    -
    \frac{2u^2-1}{u^4}  
    =
    \frac{(1-u)^4}{u^{12}}P(u),
\end{aligned}
\]
where
\[
    P(u)
    =
    u^8+6u^6-4u^5+9u^4-8u^3+6u^2-4u+1.
\]

For $u\ge1/\sqrt2$,
\[
\begin{aligned}
    P'(u)
    &=
    4\left(
        2u^7+9u^5-5u^4+9u^3-6u^2+3u-1
    \right)\\
    &=
    4\left(
        2u^7
        +u^4(9u-5)
        +3u^2(3u-2)
        +(3u-1)
    \right)>0.
\end{aligned}
\]
Hence $P$ is strictly increasing on
$[1/\sqrt2,1]$. Moreover,
\[
    P\!\left(\frac1{\sqrt2}\right)
    =
    \left(\frac94-\sqrt2\right)^2>0.
\]
Thus
\[
    P(u)>0,
    \qquad
    \frac1{\sqrt2}\le u\le1.
\]

Since $u<1$ at every finite inverse temperature,
\[
    \bigl(r_{\beta,2}^{\mathrm{Ber}}\bigr)^8>m_\beta^8.
\]
Both quantities are nonnegative, so the strict monotonicity of
$x\mapsto x^8$ on $[0,\infty)$ gives
\[
    m_\beta<r_{\beta,2}^{\mathrm{Ber}}<1.
\]
\end{proof}

Combining Proposition~\ref{prop:quench-threshold-vs-spontaneous}
with Theorem~\ref{thm:two-way-general-d} shows that, throughout the
low-temperature regime, every nonzero level
\[
    C_r,
    \qquad
    0<|r|<r_{\beta,2}^{\mathrm{Ber}},
\]
admits pathwise forcing witnesses in both directions. Since
\[
    m_\beta<r_{\beta,2}^{\mathrm{Ber}},
\]
this includes every nonzero level in the coexistence interval,
including the boundary levels \(r=\pm m_\beta\), and extends beyond it.

In particular, for
\[
    m_\beta<r<r_{\beta,2}^{\mathrm{Ber}},
\]
the biased Bernoulli quench has positive initial magnetization drift,
although its initial magnetization already exceeds the spontaneous
magnetization of the plus phase. Thus the equilibrium magnetization is
not a short-time restoring threshold for nonequilibrium biased
quenches.

\begin{remark}[Bernoulli-quench nonmonotonicity] \label {rmk:quench-monotonicity}
In dimension two, the Bernoulli-quench threshold satisfies
\[
    \beta^{\mathrm{Ber}}_2<\beta_c.
\]
Hence, if
\[
    \beta^{\mathrm{Ber}}_2<\beta\leq\beta_c,
    \qquad
    0<r<r^{\mathrm{Ber}}_{\beta,2},
\]
the Bernoulli-quench magnetization initially increases, whereas
attractiveness and uniqueness imply convergence to zero. Its magnetization curve is
therefore nonmonotone.

If $\beta>\beta_c$, Proposition~\ref{prop:quench-threshold-vs-spontaneous}
gives
\[
    m_\beta<r^{\mathrm{Ber}}_{\beta,2}.
\]
Thus, for
\[
    m_\beta<r<r^{\mathrm{Ber}}_{\beta,2},
\]
the magnetization again initially increases. By attractiveness, it is
bounded above by the magnetization of the process started from the
all-plus configuration, which converges to $m_\beta$
\cite{Liggett85,Martinelli99}. Hence its magnetization curve is nonmonotone.

Consequently, nonmonotone Bernoulli-quench magnetization occurs for
every $\beta>\beta^{\mathrm{Ber}}_2$. In the coexistence regime, the
initial magnetization may moreover be chosen above $m_\beta$.
\end{remark}

\section{The centered level and the global picture}
\label{sec:centered-level}

The slab and quench constructions both degenerate at the centered
level. For the slab family,
\[
    m'_\rho(0)
    =
    \rho\bigl(\tanh(2d\beta)-1\bigr),
\]
so the available drift vanishes as \(\rho\to0\), and no uniform drift
bound persists down to \(r=0\). In particular, the equal-width slab
has zero drift by spin-flip antisymmetry. For biased Bernoulli
quenches, the drift \(D_{\beta,d}^{\mathrm{Ber}}\) is odd, and therefore
\[
    D_{\beta,d}^{\mathrm{Ber}}(0)=0.
\]
We therefore use a separate balanced periodic configuration.

Let \(\zeta_0\) be the
\((3\mathbb Z)\times(4\mathbb Z)\)-periodic configuration on
\(\mathbb Z^2\) defined on one period cell by
\[
\begin{matrix}
+&+&-&-\\
+&-&+&-\\
+&-&-&+
\end{matrix}.
\]
The cell contains six plus spins and six minus spins. For every
\(d\geq2\), define its cylindrical extension
\(\xi_0^{(d)}\in\{-1,+1\}^{\mathbb Z^d}\) by
\[
    \xi_0^{(d)}(x_1,\ldots,x_d)
    =
    \zeta_0(x_1,x_2).
\]
It is periodic under the finite-index subgroup
\[
    \Gamma_d
    =
    (3\mathbb Z)\times(4\mathbb Z)\times\mathbb Z^{d-2},
\]
and
\[
    M_n\bigl(\xi_0^{(d)}\bigr)\longrightarrow0.
\]
Thus
\[
    \xi_0^{(d)}\in C_0.
\]

For \(t\geq0\), define
\[
    m_0^{(d)}(t)
    :=
    \frac1{12}
    \sum_{x\in\{0,1,2\}\times\{0,1,2,3\}}
    \mathbb E_{\xi_0^{(d)}}
    \bigl[\sigma_t(x_1,x_2,0,\ldots,0)\bigr].
\]
The law \(\delta_{\xi_0^{(d)}}P_t\) is invariant and ergodic under
\(\Gamma_d\). Hence the spatial ergodic theorem gives
\begin{equation}
\mathbb    P_{\xi_0^{(d)}}
    \left(
        M_n(\sigma_t)\longrightarrow m_0^{(d)}(t)
    \right)
    =1,
    \qquad t\geq0.
    \label{eq:centered-ergodic-limit}
\end{equation}
In particular,
\[
    m_0^{(d)}(0)=0.
\]

\begin{proposition}[Pathwise two-way forcing at the centered level]
\label{prop:centered-forcing}
For every \(d\geq2\) and every \(\beta>0\),
\[
    \bigl(m_0^{(d)}\bigr)'(0)\neq0.
\]
More precisely,
\[
    \bigl(m_0^{(2)}\bigr)'(0)>0,
\]
whereas
\[
    \bigl(m_0^{(d)}\bigr)'(0)<0,
    \qquad d\geq3.
\]
Consequently, \(C_0\) contains both an upward and a downward pathwise
forcing witness in every dimension \(d\geq2\).
\end{proposition}

\begin{proof}
We first record the local environments of $\zeta_0$ in dimension two.
Over one period cell, their distribution is
\[
\begin{array}{c|c|c}
\text{spin} & \text{two-dimensional local field}
            & \text{number of sites}\\
\hline
-1 & 0  & 4\\
+1 & 2  & 2\\
+1 & -2 & 2\\
-1 & 2  & 2\\
+1 & 0  & 1\\
+1 & -4 & 1.
\end{array}
\]

Set
\[
    k=d-2,
    \qquad
    T_j=\tanh(2j\beta),
    \quad j\in\mathbb Z.
\]
In the cylindrical extension, each site has $2k$ additional
neighbours, all carrying the same spin as the central site. Thus a site
with spin $s\in\{-1,+1\}$ and two-dimensional local field $h$ has
$d$-dimensional local field
\[
    h+2ks.
\]

Since the average spin over the period cell is zero, the
$-\sigma(0)$ term in the generator has zero spatial average. Averaging
the hyperbolic-tangent term over the table gives
\[
\begin{aligned}
    12\bigl(m_0^{(d)}\bigr)'(0)
    ={}&
    -4T_k+2T_{k+1}+2T_{k-1}
    +2T_{1-k}+T_k+T_{k-2}.
\end{aligned}
\]
Since $T_{-j}=-T_j$, the middle two terms cancel, and hence
\[
    12\bigl(m_0^{(d)}\bigr)'(0)
    =
    2T_{k+1}-3T_k+T_{k-2}.
\]

Define
\[
    \Delta_j:=T_{j+1}-T_j,
    \qquad j\in\mathbb Z.
\]
Then
\[
    12\bigl(m_0^{(d)}\bigr)'(0)
    =
    2\Delta_k-\Delta_{k-1}-\Delta_{k-2}.
\]
The oddness of $T_j$ gives the reflection identity
\[
    \Delta_{-1-j}=\Delta_j,
    \qquad j\in\mathbb Z.
\]
Moreover, the strict concavity of
$x\mapsto\tanh(2\beta x)$ on $[0,\infty)$ implies
\[
    \Delta_0>\Delta_1>\Delta_2>\cdots.
\]

If $k=0$, corresponding to $d=2$, the reflection identity gives
\[
\begin{aligned}
    12\bigl(m_0^{(2)}\bigr)'(0)
    =
    2\Delta_0-\Delta_{-1}-\Delta_{-2}
    =
    \Delta_0-\Delta_1>0.
\end{aligned}
\]

If $k=1$, corresponding to $d=3$, then
\[
    12\bigl(m_0^{(3)}\bigr)'(0)
    =
    2\Delta_1-\Delta_0-\Delta_{-1}
    =
    2(\Delta_1-\Delta_0)<0.
\]

If $k\ge2$, then
\[
    \Delta_{k-2}>\Delta_{k-1}>\Delta_k,
\]
and therefore
\[
    2\Delta_k-\Delta_{k-1}-\Delta_{k-2}<0.
\]

Consequently,
\[
    \bigl(m_0^{(2)}\bigr)'(0)>0,
    \qquad
    \bigl(m_0^{(d)}\bigr)'(0)<0,
    \quad d\ge3.
\]

Finally, since each $Lf_x$ is bounded and local,
$(m_0^{(d)})'$ is continuous at zero. The nonzero initial drift
therefore gives a signed linear bound on a short interval.
Equation~\eqref{eq:centered-ergodic-limit} and
Lemma~\ref{lem:pathwise-thickening} yield pathwise forcing in the
corresponding direction. Global spin-flip gives the opposite
direction.
\end{proof}

The preceding constructions yield the following theorem.

\begin{theorem}[Pathwise forcing of every exact level]
\label{thm:global-forcing}
For every \(d\ge2\), every \(\beta>0\), and every
\(r\in[-1,1]\), the level \(C_r\) contains a pathwise forcing witness.
For \(r\neq0\), forcing can be directed toward zero, while at \(r=0\)
both directions occur.

Moreover, for every \(d\ge2\) and every \(\beta>0\), there exists
\(\delta_{\beta,d}>0\) such that both forcing directions occur at every
level
\[
    r\in(-\delta_{\beta,d},\delta_{\beta,d}).
\]

If \(\beta>\beta_d^{\mathrm{Ber}}\), the outward direction can be realized
by the biased Bernoulli product measure for every
\[
    0<|r|<r_{\beta,d}^{\mathrm{Ber}}.
\]
In dimension two, throughout the phase-coexistence regime,
\[
    r_{\beta,2}^{\mathrm{Ber}}>m_\beta,
\]
so this two-way forcing interval contains
\([-m_\beta,m_\beta]\).
\end{theorem}

\begin{proof}
For every nonzero \(r\in(-1,1)\), forcing toward zero follows from
Theorem~\ref{thm:typical-toward-zero}. At nonzero rational levels,
Corollary~\ref{cor:pathwise-rational-forcing} additionally provides
an explicit periodic witness. The extremal levels \(r=\pm1\) were
treated in Section~\ref{extreme}, and the centered level follows from
Proposition~\ref{prop:centered-forcing}.

We next prove the second assertion. Choose between the centered
periodic configuration constructed above and its global spin-flip a
configuration \(\widehat\xi_0^{(d)}\) whose deterministic
magnetization curve \(\widehat m_0\) satisfies
\[
    \widehat m_0(0)=0,
    \qquad
    \widehat m_0'(0)>0.
\]
By continuity of \(\widehat m_0'\), there exist \(c,\varepsilon>0\)
such that
\[
    \widehat m_0'(t)\ge c,
    \qquad 0\le t\le\varepsilon.
\]
Set
\[
    \delta_{\beta,d}
    :=
    \widehat m_0(\varepsilon/2)>0.
\]
For every \(r\in(0,\delta_{\beta,d})\), there exists a unique
\(s_r\in(0,\varepsilon/2)\) such that
\[
    \widehat m_0(s_r)=r.
\]
Define
\[
    \nu_r^+
    :=
    \delta_{\widehat\xi_0^{(d)}}P_{s_r}
\]
The period-subgroup ergodicity argument gives
\[
    \nu_r^+(C_r)=1
\]
and, for every \(t\ge0\),
\[
    \nu_r^+P_t
    \bigl(C_{\widehat m_0(s_r+t)}\bigr)=1.
\]
Moreover, for \(0\le t\le\varepsilon/2\), we have
\(s_r+t\le\varepsilon\), and hence
\[
    \widehat m_0(s_r+t)
    \ge
    \widehat m_0(s_r)+ct
    =
    r+ct.
\]
The upward part of Corollary~\ref{forcing},
applied to \(\mu=\nu_r^+\), therefore yields a
\(\nu_r^+\)-full Borel subset of \(C_r\) consisting of upward
pathwise forcing witnesses.

Downward forcing follows from
Theorem~\ref{thm:typical-toward-zero}. Global spin-flip gives the
corresponding result for negative \(r\), while the case \(r=0\)
follows from Proposition~\ref{prop:centered-forcing}. Consequently,
for every \(\beta>0\), both forcing directions occur at every level
\[
    r\in(-\delta_{\beta,d},\delta_{\beta,d}),
\]
which proves the second assertion.

Assume now that \(\beta>\beta_d^{\mathrm{Ber}}\). For every
\(0<|r|<r_{\beta,d}^{\mathrm{Ber}}\), the Bernoulli-quench realization
of two-way forcing follows from
Theorem~\ref{thm:two-way-general-d}, while the case \(r=0\) again
follows from Proposition~\ref{prop:centered-forcing}. This proves the
third assertion. Finally, the two-dimensional statement follows from
Proposition~\ref{prop:quench-threshold-vs-spontaneous}.
\end{proof}

\section{Conclusion}

For $d\geq2$, exact magnetization levels are not dynamically homogeneous
under heat-bath Glauber dynamics. Every level $C_r$ contains a pathwise
forcing witness. At every positive temperature, two-way forcing occurs
throughout a nontrivial interval around zero.

Above the quench threshold, the outward direction is realized by the
biased Bernoulli initialization on an explicit interval. In dimension
two, throughout the coexistence regime, this interval contains
$[-m_\beta,m_\beta]$ and extends beyond the spontaneous magnetization.
In dimension two, nonmonotone Bernoulli-quench relaxation occurs for every
$\beta>\beta^{\mathrm{Ber}}_2$: for
$\beta^{\mathrm{Ber}}_2<\beta\leq\beta_c$, the magnetization initially
moves away from zero before converging back to zero, while for
$\beta>\beta_c$ one may choose the initial magnetization above
$m_\beta$ and still obtain positive initial drift.

This contrasts with dimension one, where the empirical magnetization
closes under the dynamics: $C_0$ is pathwise absorbing, while every
nonzero exact level is pathwise forced toward zero.

\bibliography{Ising}

@article{Aizenman80,
  author  = {Michael Aizenman},
  title   = {Translation Invariance and Instability of Phase Coexistence in the Two-Dimensional {I}sing System},
  journal = {Communications in Mathematical Physics},
  volume  = {73},
  number  = {1},
  pages   = {83--94},
  year    = {1980}
}

@article{Bray1994,
  author  = {Bray, A. J.},
  title   = {Theory of Phase-Ordering Kinetics},
  journal = {Advances in Physics},
  volume  = {43},
  number  = {3},
  pages   = {357--459},
  year    = {1994}
}

@article{FontesSchonmannSidoravicius02,
  author  = {Fontes, Luiz Renato and Schonmann, Roberto H. and
             Sidoravicius, Vladas},
  title   = {Stretched Exponential Fixation in Stochastic {I}sing Models
             at Zero Temperature},
  journal = {Communications in Mathematical Physics},
  volume  = {228},
  number  = {3},
  pages   = {495--518},
  year    = {2002},
  doi     = {10.1007/s002200200658}
}

@book{FriedliVelenik17,
  author    = {Sacha Friedli and Yvan Velenik},
  title     = {Statistical Mechanics of Lattice Systems: A Concrete Mathematical Introduction},
  publisher = {{C}ambridge {U}niversity {P}ress},
  address   = {Cambridge},
  year      = {2017}
}

@book{Georgii11,
  author    = {Hans-Otto Georgii},
  title     = {Gibbs Measures and Phase Transitions},
  edition   = {2},
  publisher = {De Gruyter},
  address   = {Berlin/Boston},
  year      = {2011},
  series    = {De Gruyter Studies in Mathematics},
  volume    = {9}
}

@article{GheissariSinclair23,
  author  = {Gheissari, Reza and Sinclair, Alistair},
  title   = {Low-temperature {Ising} dynamics with random initializations},
  journal = {The Annals of Applied Probability},
  volume  = {33},
  number  = {5},
  pages   = {3916--3957},
  year    = {2023},
  doi     = {10.1214/22-AAP1911}
}

@article{GheissariSly26,
  author  = {Gheissari, Reza and Sly, Allan},
  title   = {Rapid Phase Ordering of {I}sing Dynamics on
             $\mathbb{Z}^2$},
  journal = {arXiv},
  volume  = {arXiv:2605.08052},
  year    = {2026}
}

@article{Glauber63,
  author  = {Glauber, Roy J.},
  title   = {Time-Dependent Statistics of the {I}sing Model},
  journal = {Journal of Mathematical Physics},
  volume  = {4},
  number  = {2},
  pages   = {294--307},
  year    = {1963}
}

@article{HaggstromSteif00,
  author  = {H{\"a}ggstr{\"o}m, Olle and Steif, Jeffrey E.},
  title   = {Propp--{W}ilson Algorithms and Finitary Codings for High Noise
             {M}arkov Random Fields},
  journal = {Combinatorics, Probability and Computing},
  volume  = {9},
  number  = {5},
  pages   = {425--439},
  year    = {2000},
  doi     = {10.1017/S0963548300004363}
}

@article{Harris72,
  author  = {Theodore E. Harris},
  title   = {Nearest-neighbor {M}arkov interaction processes on multidimensional lattices},
  journal = {Advances in Mathematics},
  volume  = {9},
  number  = {1},
  pages   = {66--89},
  year    = {1972}
}

@incollection{Higuchi81,
  author    = {Yasunari Higuchi},
  title     = {On the Absence of Non-Translation Invariant {G}ibbs States for the Two-Dimensional {I}sing System},
  booktitle = {Random Fields, Vol. II},
  series    = {Colloquia Mathematica Societatis J{\'a}nos Bolyai},
  volume    = {27},
  publisher = {North-Holland},
  address   = {Amsterdam},
  pages     = {517--534},
  year      = {1981}
}

@article{LeNyRedig02,
  author  = {Le Ny, Arnaud and Redig, Frank},
  title   = {Short-Time Conservation of Gibbsianness under Local
             Stochastic Evolutions},
  journal = {Journal of Statistical Physics},
  volume  = {109},
  number  = {5--6},
  pages   = {1073--1090},
  year    = {2002},
  doi     = {10.1023/A:1020418109483}
}

@book{Liggett85,
  author    = {Liggett, Thomas M.},
  title     = {Interacting Particle Systems},
  series    = {Grundlehren der Mathematischen Wissenschaften},
  volume    = {276},
  publisher = {Springer-Verlag},
  address   = {New York},
  year      = {1985}
}

@book{Liggett99,
  author    = {Liggett, Thomas M.},
  title     = {Stochastic Interacting Systems: Contact, Voter and Exclusion Processes},
  series    = {Grundlehren der Mathematischen Wissenschaften},
  volume    = {324},
  publisher = {Springer-Verlag},
  address   = {Berlin},
  year      = {1999}
}

@book{Lothaire02,
  author    = {Lothaire, Monsieur},
  title     = {Algebraic Combinatorics on Words},
  series    = {Encyclopedia of Mathematics and its Applications},
  volume    = {90},
  publisher = {{C}ambridge {U}niversity {P}ress},
  address   = {{C}ambridge},
  year      = {2002}
}

@article{LubetzkySly16,
  author  = {Lubetzky, Eyal and Sly, Allan},
  title   = {Information Percolation and Cutoff for the Stochastic {I}sing Model},
  journal = {Journal of the American Mathematical Society},
  volume  = {29},
  number  = {3},
  pages   = {729--774},
  year    = {2016}
}

@incollection{Martinelli99,
  author    = {Fabio Martinelli},
  title     = {Lectures on {G}lauber Dynamics for Discrete Spin Models},
  booktitle = {Lectures on Probability Theory and Statistics},
  series    = {Lecture Notes in Mathematics},
  volume    = {1717},
  publisher = {Springer},
  address   = {Berlin, Heidelberg},
  year      = {1999},
  pages     = {93--191},
  doi       = {10.1007/BFb0092673}
}

@article{vanEnterFernandezdenHollanderRedig02,
  author  = {van Enter, Aernout C. D. and Fern{\'a}ndez, Roberto and
             den Hollander, Frank and Redig, Frank},
  title   = {Possible Loss and Recovery of Gibbsianness during the
             Stochastic Evolution of {G}ibbs Measures},
  journal = {Communications in Mathematical Physics},
  volume  = {226},
  number  = {1},
  pages   = {101--130},
  year    = {2002},
  doi     = {10.1007/s002200200605}
}

@book{Walters82,
  author    = {Walters, Peter},
  title     = {An Introduction to Ergodic Theory},
  series    = {Graduate Texts in Mathematics},
  volume    = {79},
  publisher = {Springer},
  address   = {New York},
  year      = {1982}
}

@article{Yang52,
  author  = {Yang, Chen-Ning},
  title   = {The Spontaneous Magnetization of a Two-Dimensional {I}sing Model},
  journal = {Physical Review},
  volume  = {85},
  number  = {5},
  pages   = {808--816},
  year    = {1952}
}

\end{document}